# ON THE LATTICE OF COTORSION THEORIES

RÜDIGER GÖBEL, SAHARON SHELAH, AND SIMONE L. WALLUTIS

ABSTRACT. We discuss the lattice of cotorsion theories for abelian groups. First we show that the sublattice of the well–studied rational cotorsion theories can be identified with the well–known lattice of types. Using a recently developed method for making Ext vanish we also prove that any power set together with the ordinary set inclusion (and thus any poset) can be embedded into the lattice of all cotorsion theories.

## INTRODUCTION

Throughout this paper we work in the category Mod–$\mathbb{Z}$ of abelian groups although most of the notions and results can be extended to module categories over arbitrary rings, in particular everything remains true for modules over principal ideal domains.

Cotorsion theories for abelian groups have been introduced by Salce in 1979 [17]. Following his notation we call a a pair $(\mathcal{F}, \mathcal{C})$ a *cotorsion theory* if $\mathcal{F}$ and $\mathcal{C}$ are classes of abelian groups which are maximal with respect to the property that $\text{Ext}(F, C) = 0$ for all $F \in \mathcal{F}$, $C \in \mathcal{C}$. Salce [17] has shown that every cotorsion theory is cogenerated by a class of groups that are torsion or torsion–free where $(\mathcal{F}, \mathcal{C})$ is said to be cogenerated by the class $\mathcal{A}$ if $\mathcal{C} = \mathcal{A}^\perp = \{X \mid \text{Ext}(A, X) = 0 \text{ for all } A \in \mathcal{A}\}$ and $\mathcal{F} = {}^\perp(\mathcal{A}^\perp) = \{Y \mid \text{Ext}(Y, X) = 0 \text{ for all } X \in \mathcal{A}^\perp\}$. Very recently Bican, El Bashir and Enochs [1] have shown that the flat cotorsion theory is actually cogenerated by a set – a result which has been unknown for the last two decades and which plays a crucial role in proving the flat cover conjecture. Note that being cogenerated by a set is the same as being cogenerated by a single group although this group is likely to be mixed as a direct sum of all groups from the cogenerating set. However, it is not known whether every cotorsion theory is singly cogenerated or not.

Other results on cotorsion theories have been proven, for example, concerning the existence of enough projectives and enough injectives. The basic work has been done by Salce in his original paper [17] where among other results he has proven that a cotorsion theory has enough injectives if and only if it has enough projectives. In a quite recent paper the first and the second author developed a method to construct splitters, that is groups $G$ satisfying $\text{Ext}(G, G) = 0$, which could be applied to prove the existence of enough projectives, respectively enough injectives in so–called rational cotorsion theories (see [14]). This method has been improved independently by Eklof, Trlifaj [6] and the last author [16] and will also be used here.

This work is supported by GIF project No. G–0545–173.06/97 of the German – Israeli Foundation for Scientific Research & Development.
It is [GPSh:721] in Shelah's list of publications;
Third author's maiden name is Pabst; supported by Deutsche Forschungsgemeinschaft.

However, in this paper we shall discuss the lattice structure of the class of all cotorsion theories. We order the pairs correspondingly to the second component $\mathcal{C}$, the so–called cotorsion class; the first component $\mathcal{F}$ is said to be the torsion–free class. We say that $(\mathcal{F}, \mathcal{C}) \leq (\mathcal{F}', \mathcal{C}')$ if $\mathcal{C} \subseteq \mathcal{C}'$ or, equivalently, $\mathcal{F} \supseteq \mathcal{F}'$. The minimal element with respect to this order is $(\mathrm{Mod}{-}\mathbb{Z}, \mathcal{D})$ where $\mathcal{D}$ is the class of all divisible groups; it is, for example, cogenerated by the set of all cyclic groups of prime order or, equivalently, by the single group $\bigoplus_{p \in \Pi} Z_p$ ($\Pi =$ "all primes"). The maximal element is the cotorsion theory $(\mathcal{L}, \mathrm{Mod}{-}\mathbb{Z})$ where $\mathcal{L}$ is the class of all free groups; it is cogenerated by $\mathbb{Z}$. Another important and well–studied cotorsion theory is the classical one $(\mathcal{F}_0, \mathcal{C}_0)$ where $\mathcal{F}_0$ denotes the class of all torsion–free groups and $\mathcal{C}_0$ denotes the class of all cotorsion groups; it is cogenerated by the rationals $\mathbb{Q}$. Canonically we define the infimum and supremum of a family $\{(\mathcal{F}_i, \mathcal{C}_i) \mid i \in I\}$ of cotorsion theories by

$$\bigwedge_{i \in I} (\mathcal{F}_i, \mathcal{C}_i) = \left( {}^\perp(\bigcap_{i \in I} \mathcal{C}_i), \bigcap_{i \in I} \mathcal{C}_i \right) \text{ and } \bigvee_{i \in I} (\mathcal{F}_i, \mathcal{C}_i) = \left( \bigcap_{i \in I} \mathcal{F}_i, (\bigcap_{i \in I} \mathcal{F}_i)^\perp \right).$$

As we have said before every cotorsion theory is cogenerated by a class of groups which are torsion or torsion–free and hence it is the infimum of a cotorsion theory cogenerated by torsion groups and a cotorsion theory cogenerated by torsion–free groups. Now, a cotorsion theory which is cogenerated by torsion groups is always less than or equal to the classical one. But the sublattice of all cotorsion theories between the minimal and the classical one has already been characterized by Salce [17, Proposition 2.8] and therefore we can restrict our attention to the cotorsion theories cogenerated by torsion–free groups, i.e. to the cotorsion theories above the classical one.

Naturally, we first consider easy cases of torsion–free groups, namely the rank–1 groups which are also called rational groups as they can be identified with the subgroups of the rationals. Corresponding to the latter we call a cotorsion theory cogenerated by a rank–1 group a *rational cotorsion theory*. Rational cotorsion theories have been discussed in detail by Salce [17]. Using his characterization we shall prove in §1 that the sublattice of all rational cotorsion theories can be identified with the well–known lattice of types. We construct examples to establish that the (obvious) lattice anti–homomorphism is an anti–isomorphism. Note that for the lattice of types it is known that there exist anti–chains of size $2^{\aleph_0}$ which equals the cardinality of the lattice (see [10]) and also ascending and descending chains of uncountable length, in fact there are descending and ascending chains of cofinality at least $\aleph_1$ (see [4]).

Turning our attention to the more general case of all cotorsion theories, as we shall do in §3, we cannot find any obvious "candidate" with which the lattice could be compared; note that there is a proper class of cotorsion theories. However, knowing about the properties of the lattice of types or, equivalently, of the lattice of rational cotorsion theories as mentioned above it seems natural to ask if there exist ascending, descending, and anti–chains of cotorsion theories of arbitrary size. The existence of descending chains of arbitrary length follows immediately from results proven by the first author and Trlifaj (see [15]). However, we can prove that there is an affirmative answer to the above question in general. Actually, to our own surprise, we can show even more:

**Main Theorem 3.1** *Any power set* $(\mathcal{P}, \subseteq)$ *can be embedded into the lattice of all cotorsion theories.*

Therefore any partial order can be embedded into the lattice of all cotorsion theories. In fact, any poset can be embedded into the lattice of all singly cogenerated cotorsion theories.

In order to prove the Main Theorem we shall construct groups $G_X$ and $H^X$ for any subset $X$ of an arbitrary but fixed set $I$ such that $\text{Ext}(G_X, H^Y) = 0$ if and only if $X \subseteq Y$. This way we obtain an order–reversing and injective morphism from $(\mathcal{P}(I), \subseteq)$ into the lattice of all cotorsion theories by mapping the set $X$ onto the cotorsion theory cogenerated by $G_X$. Of course, this implies the required embedding since $(\mathcal{P}, \subseteq)$ is anti–isomorphic to itself.

It is known how to construct a group $H$ such that $\text{Ext}(G, H) = 0$ for a given group $G$ or even for a collection of groups. As mentioned before this method has been introduced in [14] and futher developed in [6] and [16]. We shall use this method as presented in [6] in §3 to construct the groups $H^X$ ($X \subseteq I$). For this construction it is not important what the groups $G_X$ look like. However, the $G_X$'s need to satisfy certain properties to guarantee that Ext is non zero in some cases. We find this is amazing as one would expect that it is obvious how to get non–vanishing Ext's after all the hard work which had been done over decades in order to establish a method for making Ext vanish. But we have also some work to do to obtain non–zero Ext's. The key to prove $\text{Ext}(G_X, H^Y) \neq 0$ for $X \nsubseteq Y$ is the existence of a stationary set $S$ such that $H^Y$ is locally $S$–free and $G_X$ is not, where a group $A$ is said to be *locally $S$–free* for a stationary set $S$ of a cardinal $\kappa$ if, for any chain $\{K_\alpha \mid \alpha < \kappa\}$ of subgroups $K_\alpha$ of $A$ with $|K_\alpha| < \kappa$ ($\alpha < \kappa$), the set $\{\delta \in S \mid K_{\delta+1}/K_\delta \text{ not } \aleph_1\text{–free}\}$ is not stationary in $\kappa$ (see Definition 3.4). Therefore we have to construct the groups $G_X$ in such a way that they are locally $S$–free with respect to some stationary set $S$ and not locally $S'$–free with respect to others. In fact, we shall construct a group $G$ depending on a stationary set $S$ and then we define the groups $G_X$ depending on different stationary sets. The construction of $G$ is interesting on its own and thus we consider it in a separate section; in §2 we will prove:

**Theorem 2.8** *Let $\kappa$ be a regular cardinal with $|\alpha|^{\aleph_0} < \kappa$ for all $\alpha < \kappa$. Then there exists an $\aleph_1$–free group $G$ of cardinality $\kappa$ such that, for any subgroup $U$ of $G$ with $|U| = \kappa$, either $G = U$ or $G/U$ is not cotorsion–free.*

The corresponding result for strong limit singular cardinals has been proven by the first and second authors ([11]) at a stage when the Black Box was not yet fully developed. They already suggested that the above result is true but it hadn't been proven so far.

To see the connection with a given stationary set $S$ let us finally note that the group $G$ which will be constructed to prove Theorem 2.8 has a $\kappa$–filtration $G = \bigcup_{\alpha < \kappa} G_\alpha$ such that $G/G_\alpha$ is $\aleph_1$–free if and only if $\alpha \in S$. Using this fact we can show that $G$ depending on $S$ is locally $S'$–free for any stationary set $S'$ disjoint from $S$ but not locally $S$–free. While the construction of $G$ will be provided in §2 we shall prove the latter in §3.

However, we begin with some "warm–ups", namely with the rational cotorsion theories.

## 1. The lattice of all rational cotorsion theories

In this section we describe the lattice of all cotorsion theories which are cogenerated by a rank–1 group. For a rank–1 group $T$ let $T^\perp = \{X \mid \text{Ext}(T, X) = 0\}$ and $^\perp(T^\perp) = \{Y \mid \text{Ext}(Y, X) = 0 \text{ for all } X \in T^\perp\}$. The pair $\left(^\perp(T^\perp),\, T^\perp\right)$ is the cotorsion theory cogenerated by $T$ where $T^\perp$ is the class of all $T$–cotorsion groups and $^\perp(T^\perp)$ the corresponding torsion–free class. Since rank–1 groups are also called rational groups, as rank–1 groups can be identified with the subgroups of the rationals $\mathbb{Q}$, we refer to the cotorsion theory $\left(^\perp(T^\perp),\, T^\perp\right)$ as a *rational cotorsion theory*. Rational cotorsion theories have been discussed by Salce [17] in detail. We shall use his results to establish an order–reversing isomorphism between the lattice of rational cotorsion theories and the lattice of types. As cotorsion theories in general we also order the rational cotorsion theories according to the inclusion of the cotorsion classes (see §0).

We think of a type $\tau$ as a sequence $(t_p)_{p\in\Pi}$, where $t_p \in \mathbb{N} \cup \{\infty\}$ and $\Pi$ is the set of all primes in $\mathbb{Z}$, keeping in mind that this sequence represents an equivalence class. Recall that two such sequences are equivalent if they only differ in finitely many finite entries.

It is well known that rank–1 groups are uniquely determined by their types (up to isomorphism) and that there exists a rank–1 group of type $\tau$ for each possible type $\tau$. For more details we refer to [10].

Now let $\mathcal{T}$ be the set of all types and let $\mathcal{C}_{rat}$ be the set of all rational cotorsion theories. We define $\Phi : (\mathcal{T}, \leq) \longrightarrow (\mathcal{C}_{rat}, \leq)$ by $\tau\Phi = \left(^\perp(T^\perp),\, T^\perp\right) \in \mathcal{C}_{rat}$ where $T$ is a rank–1 group of type $\tau$. The aim of this section is to prove that the mapping $\Phi$ is an order–reversing isomorphism.

First we show that $\Phi$ is order reversing.

**Lemma 1.1.** *Let $T$, $R$ be rank–1 groups with $t(T) \leq t(R)$. Then $T^\perp \supseteq R^\perp$.*

*Proof.* Since $t(T) \leq t(R)$ there is a monomorphism $\varepsilon : T \longrightarrow R$ (see [10, Proposition 85.4]). Let $G$ be an element of $R^\perp$, i.e. $\mathrm{Ext}(R,G) = 0$. Now, the short exact sequence
$$0 \longrightarrow T \xrightarrow{\varepsilon} R \longrightarrow R/T\varepsilon \longrightarrow 0$$
induces the exact sequence
$$\mathrm{Ext}(R/T\varepsilon, G) \longrightarrow \mathrm{Ext}(R, G) \longrightarrow \mathrm{Ext}(T, G) \longrightarrow 0$$
and hence $\mathrm{Ext}(T, G) = 0$. Therefore $G \in T^\perp$ and thus $R^\perp \subseteq T^\perp$. □

It follows immediately from the above lemma that the mapping $\Phi$ is well defined:

**Corollary 1.2.** *Let $T$, $T'$ be rank–1 groups of the same type $\tau$. Then the corresponding cotorsion classes $T^\perp$ and $(T')^\perp$ coincide.*

Note, that more generally we have $G \subseteq H$ implies $G^\perp \supseteq H^\perp$ and $G \cong H$ implies $G^\perp = H^\perp$ for any groups $G$, $H$.

In order to show that $\Phi$ is an isomorphism we consider types $\tau = t(T)$, $\rho = t(R)$ with $\tau$ strictly less than $\rho$ and we show that $R^\perp$ is properly contained in $T^\perp$, i.e. we construct groups $G \in T^\perp \setminus R^\perp$.

Throughout the remainder of this section let $\tau = (t_p)_{p\in\Pi} = t(T)$ and $\rho = (r_p)_{p\in\Pi} = t(R)$ with $t_p \leq r_p$ for all primes $p$. For $\tau$ to be strictly less than $\rho$ one of the following two conditions has to be satisfied:

(1) There exists a prime $q$ such that $t_q < \infty$ and $r_q = \infty$.
(2) There is an infinite set $P$ of primes such that $t_p < r_p < \infty$ for all $p \in P$.

Before we can construct the required groups we need some properties of $T$–cotorsion groups. Fortunately, $T$–cotorsion groups have already been characterized by Salce [17, Theorem 3.5]:

**Proposition 1.3.** *Let $\tau = t(T)$ be as above. Then*
$$G \in T^\perp \iff G/G_\tau \cong \prod_{p \in \Pi} G_p^\tau \iff G/G_\tau \text{ is } (\mathbb{Q}\text{-})\text{cotorsion}$$
*where* $G_\tau = \bigcap_{p \in \Pi} p^{t_p} G$, $G_p^\tau = G/p^{t_p} G$ *for* $t_p < \infty$, *and*
$G_p^\tau = \text{Ext}(Z_{p^\infty}, G)$ *for* $t_p = \infty$.

Applying the above proposition to rank–1 groups gives the following:

**Corollary 1.4.** *Let $X$ be a rational group with $t(X) = (x_p)_{p \in \Pi}$ and let $\tau = t(T)$ be as above.*
*Then $X$ is an element of $T^\perp$ if and only if $x_p = \infty$ for almost all $p$ with $t_p \neq 0$ and whenever $t_p = \infty$.*

*Proof.* First recall that, for an abelian group $G$, $\text{Ext}(Z_{p^n}, G) \cong G/p^n G$ for any $n \in \mathbb{N}$ and $\text{Ext}(Z_{p^\infty}, G) \cong \bigoplus_m J_p$ where $m$ is the rank of a $p$–basic subgroup of $G$ and $J_p$ is the additive group of the ring of $p$–adic integers.

Now assume that $X \in T^\perp$. Then $X/X_\tau \cong \prod_{\{p \in \Pi \mid t_p \neq 0\}} \text{Ext}(Z_{p^{t_p}}, X)$ by Proposition 1.3. But $X$ and hence $X/X_\tau$ is countable and thus $\text{Ext}(Z_{p^{t_p}}, X) = 0$ for almost all $p$ with $t_p \neq 0$ and whenever $t_p = \infty$ since $\left|\prod_{n \in \omega} M_n\right| \geq 2^{\aleph_0}$ if $|M_n| \geq 2$ for all $n$, and $J_p \subseteq \text{Ext}(Z_{p^\infty}, X)$ unless $\text{Ext}(Z_{p^\infty}, X) = 0$. Therefore $p^{t_p} X = X$ for almost all $p$ with $0 < t_p < \infty$ and the rank of a $p$–basic subgroup of $X$ is zero for all $p$ with $t_p = \infty$. In either case it follows that $X$ is $p$–divisible and hence $x_p = \infty$ for almost all $p$ with $t_p \neq 0$ and whenever $t_p = \infty$. □

From the above it is clear that the primes with $t_p = 0$ play a special role. In particular, it makes sense to divide the second case (2) into two subcases:

**(2a)** There is an infinite set $P$ of primes such that $t_p = 0$ and $0 \neq r_p < \infty$ for all $p \in P$.

**(2b)** There is an infinite set $P$ of primes such that $0 < t_p < r_p < \infty$ for all $p \in P$.

However, we first consider case (1).

**Proposition 1.5.** *Suppose $t(T) = \tau < \rho = t(R)$ such that (1) is satisfied. Then there exists a rank–1 group $X$ which is an element of $T^\perp$ but not of $R^\perp$.*

*Proof.* Suppose (1), i.e. there is a prime $q$ such that $t_q < \infty$ and $r_q = \infty$.
Let $X = \mathbb{Z}_{(q)} = \{\frac{m}{n} \in \mathbb{Q} \mid (n, q) = 1\}$ be the localization of the integers $\mathbb{Z}$ at the prime $q$, i.e. $t(X) = (x_p)_{p \in \Pi}$ with $x_p = \infty$ for all $p \neq q$ and $x_q = 0$. Then $X \in T^\perp \setminus R^\perp$ by Corollary 1.4. □

Case (2a) is as easily tackled as the above:

**Proposition 1.6.** *Suppose $t(T) = \tau < \rho = t(R)$ such that (2a) is satisfied. Then there exists a rank–1 group $X$ which is an element of $T^\perp$ but not of $R^\perp$.*

*Proof.* Suppose (2a), i.e. there is an infinite set $P$ of primes $p$ with $t_p = 0$ and $0 \neq r_p < \infty$.
We define $X \subseteq \mathbb{Q}$ by $t(X) = (x_p)_{p \in \Pi}$ with $x_p = 1$ for $p \in P$ and $x_p = \infty$ otherwise, i.e $X = \left\langle \{\frac{1}{p} \mid p \in P\} \cup \{\frac{1}{p^n} \mid n \in \omega, p \in \Pi \setminus P\} \right\rangle$. Then $G$ is an element of $T^\perp \setminus R^\perp$ by Corollary 1.4. □

It remains to consider the case (2b). This is slightly more difficult as we cannot expect to find a rank–1 group belonging to $T^\perp$ but not to $R^\perp$ by Corollary 1.4. In fact, we cannot even find a group of any finite rank belonging to $T^\perp$ and not to $R^\perp$ as we shall see shortly. Beforehand we need:

**Lemma 1.7.** *Let $T$ and $X$ be rank–1 groups with $\operatorname{Ext}(T,X) \neq 0$. Then $\operatorname{Ext}(T,X)$ has cardinality $2^{\aleph_0}$.*

*Proof.* Let $t(T) = \tau = (t_p)_{p \in \Pi}$ and $t(X) = (x_p)_{p \in \Pi}$ as before. The short exact sequence
$$0 \longrightarrow \mathbb{Z} \longrightarrow T \longrightarrow T/\mathbb{Z} \longrightarrow 0$$
induces the exact sequence
$$\operatorname{Hom}(\mathbb{Z},X) \longrightarrow \operatorname{Ext}(T/\mathbb{Z},X) \longrightarrow \operatorname{Ext}(T,X) \longrightarrow \operatorname{Ext}(\mathbb{Z},X) = 0.$$
Now, $T/\mathbb{Z} \cong \bigoplus_{t_p \neq 0} Z_{p^{t_p}}$ and hence $\operatorname{Ext}(T/\mathbb{Z},X) \cong \prod_{t_p \neq 0} \operatorname{Ext}(Z_{p^{t_p}}, X) =: E$.

By assumption $\operatorname{Ext}(T,X) \neq 0$ and thus, by Corollary 1.4, there is either some prime $q$ with $x_q < \infty$ and $t_q = \infty$ or there are infinitely many primes $q_n$ with $x_{q_n} < \infty$ and $0 \neq t_{q_n} < \infty$ $(n \in \omega)$. In the first case, $\operatorname{Ext}(Z_{q^\infty}, X)$ contains a copy of the $q$–adic integers $J_q$ and in the latter we have $\operatorname{Ext}\left(Z_{q_n^{t_{q_n}}}, X\right) \cong X\left/q_n^{t_{q_n}} X\right. \neq 0$ for all $n \in \omega$. Therefore, in either case it follows that $E$ has at least cardinality $2^{\aleph_0}$ and, of course, the cardinality cannot be bigger.

Finally, the mapping $\operatorname{Ext}(T/\mathbb{Z},X) \longrightarrow \operatorname{Ext}(T,X)$ in the above sequence is an epimorphism with at most countable kernel and thus the result follows. □

Now we can proceed with:

**Proposition 1.8.** *Let $t(T) = \tau < \rho = t(R)$ satisfying condition (2b) but neither (1) nor (2a) and let $\mathcal{F}$ denote the set of all finite rank torsion–free groups. Then $T^\perp \cap \mathcal{F} = R^\perp \cap \mathcal{F}$.*

*Proof.* Without loss of generality we may assume that $t_p = 0$ iff $r_p = 0$, $t_p = \infty$ iff $r_p = \infty$ and that the remaining set is $P = \{p \in \Pi \mid 0 < t_p < r_p < \infty\}$ which is infinite by assumption. Obviously $T^\perp \cap \mathcal{F} \supseteq R^\perp \cap \mathcal{F}$. So let $G \in T^\perp \cap \mathcal{F}$ be of rank $n$. We show $G \in R^\perp$ by induction on $n$.

For $n = 1$ this follows immediately from Corollary 1.4. So, let $n > 1$ and consider the short exact sequence
$$0 \longrightarrow X \longrightarrow G \longrightarrow G/X \longrightarrow 0$$
where $X$ is a pure subgroup of $G$ of rank 1 and so $G/X$ is torsion–free of rank $n-1$. The above sequence induces the exact sequences
$$\operatorname{Hom}(T, G/X) \longrightarrow \operatorname{Ext}(T,X) \longrightarrow \operatorname{Ext}(T,G) \longrightarrow \operatorname{Ext}(T,G/X) \longrightarrow 0$$
and
$$\operatorname{Hom}(R, G/X) \longrightarrow \operatorname{Ext}(R,X) \longrightarrow \operatorname{Ext}(R,G) \longrightarrow \operatorname{Ext}(R,G/X) \longrightarrow 0.$$
Now $\operatorname{Ext}(T,G) = 0$ by assumption and so also $\operatorname{Ext}(T,G/X) = 0$. Hence, by induction hypothesis, $\operatorname{Ext}(R,G/X) = 0$ since $rk(G/X) = n-1$. Henceforce, the first of the two above sequences reduces to
$$\operatorname{Hom}(T,G/X) \longrightarrow \operatorname{Ext}(T,X) \longrightarrow 0$$
and so $|\operatorname{Ext}(T,X)| \leq \aleph_0$. But this is only possible if $\operatorname{Ext}(T,X) = 0$ by Lemma 1.7. Therefore we also have $\operatorname{Ext}(R,X) = 0$ since $rk(X) = 1$ and thus it follows from the second of the above sequences that $\operatorname{Ext}(R,G) = 0$, i.e. $G \in R^\perp$. □

The above proposition shows that in order to get a counterexample for case (2b) the required group $G$ needs to be of infinite rank. We would like to thank Luigi Salce for suggesting this.

Before we tackle the remaining case (2b) we recall some well–known facts:

**Remark 1.9.** (a) *Algebraically compact groups, in particular complete groups, are cotorsion.*
  (b) *Let $H$ be a reduced cotorsion group and let $U \subseteq H$ be a subgroup. Then $U$ is cotorsion if and only if $H/U$ is reduced.*
  (c) *The completion of $\bigoplus_{p \in P} Z_{p^{n_p}}$ in the $\mathbb{Z}$–adic topology is $\prod_{p \in P} Z_{p^{n_p}}$ for any set $P$ of primes.*

Part (c) is an easy exercise and is left to the reader; (a), (b) can be found in [9, pp. 232/233].

Now we are ready for:

**Proposition 1.10.** *Suppose $t(T) = \tau < \rho = t(R)$ such that (2b) is satisfied. Then there exists a group $G$ which is an element of $T^\perp$ but not of $R^\perp$.*

*Proof.* Suppose (2b), i.e. there exists an infinite set $P$ of primes such that $0 < t_p < r_p < \infty$.

Let $H = \bigoplus_{p \in P} \mathbb{Z}_{(p)} \subseteq_* \prod_{p \in P} \mathbb{Z}_{(p)} = H'$ where $\mathbb{Z}_{(p)}$ is the localization of the integers at the prime $p$. Note that $H'$ and thus also its pure subgroup $H$ is $q$–divisible for any prime $q \notin P$. We define $G$ as a subset of $H'$ by $G =$
$\{(g_p)_{p \in P} \in H' \mid \exists m, k \in \mathbb{N} \text{ s.t. } mg_p \in \mathbb{Z} \text{ and } |mg_p| \leq kp^{t_p} \text{ for all } p \in P\}.$

First we show that $G$ is a pure subgroup of $H'$ containing $H$.

Let $(g_p)_{p \in P}, (h_p)_{p \in P} \in G$, i.e. there are $m, n, k, l \in \mathbb{N}$ such that $mg_p, nh_p \in \mathbb{Z}$, $|mg_p| \leq kp^{t_p}$, and $|nh_p| \leq lp^{t_p}$ for all $p \in P$. Then $mn(g_p + h_p) = n(mg_p) + m(nh_p) \in \mathbb{Z}$ and $|mn(g_p + h_p)| \leq n|mg_p| + m|nh_p| \leq nkp^{t_p} + mlp^{t_p} = (nk + ml)p^{t_p}$ for all $p \in P$. Thus $(g_p)_{p \in P} + (h_p)_{p \in P} \in G$, i.e. $G$ is a subgroup. As an immediate consequence from the definition we have that $G$ is pure in $H'$.

Now let $(h_p)_{p \in P}$ be an element of $H$, i.e. $h_p = 0$ for almost all $p$ and $h_p = \dfrac{z_p}{n_p} \in \mathbb{Z}_{(p)}$. Let $N$ be the product of all $n_p$ and let $K$ be the sum of all $|h_p|$ over all $p$ with $h_p \neq 0$. Then $Nh_p \in \mathbb{Z}$ and $N|h_p| \leq NK \; (\in \mathbb{Z})$ for all $p \in P$. Thus $H \subseteq G$.

Next let $\pi_t : H' \longrightarrow H'_t = \prod_{p \in P} \mathbb{Z}_{(p)} / p^{t_p} \mathbb{Z}_{(p)}$ be the canonical epimorphism given by $(h_p)_{p \in P} \pi_t = \left(h_p + p^{t_p} \mathbb{Z}_{(p)}\right)_{p \in P}$. Obviously, $H\pi_t \cong \bigoplus_{p \in P} \mathbb{Z}_{(p)} / p^{t_p} \mathbb{Z}_{(p)} \cong \bigoplus_{p \in P} Z_{p^{t_p}}$ and thus $H'_t \cong \prod_{p \in P} Z_{p^{t_p}}$ is the completion of $H\pi_t$ by 1.9 (c). Therefore $H'_t$ is cotorsion by 1.9 (a).

Since an element of $H'_t$ can be represented by $(g_p)_{p \in P}$ with $g_p \in \mathbb{Z}$ and $0 \leq g_p < p^{t_p}$ we also have immediately that $G\pi_t = H'_t$. But $G\pi_t \cong G \Big/ \bigcap_{p \in \Pi} p^{t_p} G$ since $G$ is a pure subgroup of $H'$ and so $\ker \pi_t \cap G =$
$\left(\prod_{p \in P} p^{t_p} \mathbb{Z}_{(p)}\right) \cap G = \left(\bigcap_{p \in P} p^{t_p} H'\right) \cap G = \bigcap_{p \in P} (p^{t_p} H' \cap G) = \bigcap_{p \in P} p^{t_p} G = \bigcap_{p \in \Pi} p^{t_p} G = G_\tau$. Thus we have shown that $G/G_\tau \cong G\pi_t = H'_t$ is cotorsion. Therefore $G$ is an element of $T^\perp$ by Proposition 1.3.

Finally we show that $G$ is not an element of $R^\perp$. Following the same arguments as above we have that $H'_r = \prod_{p \in P} \mathbb{Z}_{(p)} / p^{r_p} \mathbb{Z}_{(p)}$ is the completion of $H\pi_r \cong \bigoplus_{p \in P} Z_{p^{r_p}}$, that $G\pi_r \cong G/G_\rho$ and that $H'_r$ is cotorsion where $G_\rho = \bigcap_{p \in \Pi} p^{r_p} G$ and $\pi_r : H' \longrightarrow H'_r$ is the corresponding epimorphism.

Now $H\pi_r \subseteq G\pi_r$ and so $H'_r/G\pi_r$ is divisible as an epimorphic image of the divisible group $H'_r/H\pi_r$. Therefore it is enough to show that $G\pi_r \neq H'_r$ to prove $G \notin R^\perp$ by 1.9 (b).

We choose integers $n_p$ $(p \in P)$ such that $p^{t_p+\frac{1}{2}} - 1 \leq n_p \leq p^{t_p+\frac{1}{2}}$. Suppose $(n_p + p^{r_p}\mathbb{Z}_{(p)})_{p \in P} \in G\pi_r$. Then there is $(g_p)_{p \in P} \in G$ such that $n_p \equiv g_p \mod p^{r_p}\mathbb{Z}_{(p)}$ for all $p \in P$. Note that $g_p \neq n_p$ for almost all $p$ since $p^{t_p+\frac{1}{2}} - 1 \leq n_p \leq kp^{t_p}$ for all $p \in P$ is impossible for a fixed $k \in \mathbb{N}$. However, there are $m, k \in \mathbb{N}$ such that $mg_p \in \mathbb{Z}$ and $|mg_p| \leq kp^{t_p}$ for all $p \in P$. So $mn_p \equiv mg_p \mod p^{r_p}\mathbb{Z}$ and thus $p^{r_p}$ divides $m(n_p - g_p)$ (in $\mathbb{Z}$). Since $t_p < r_p$ we have that $p^{t_p+1}$ divides $|m(n_p - g_p)| \leq |mn_p| + |mg_p| \leq mp^{t_p+\frac{1}{2}} + mkp^{t_p} = mp^{t_p}(p^{\frac{1}{2}} + k)$ for all $p \in P$. But, for almost all $p \in P$, $m, k < \frac{1}{2}p^{\frac{1}{2}}$ and thus $p^{t_p}(mp^{\frac{1}{2}} + mk) < p^{t_p+1}$ contradicting $m(n_p - g_p) \neq 0$. Therefore $(n_p)_{p \in P} \in H'_r \setminus G\pi_r$ and this completes the proof. □

As a consequence of the above results we can finally state:

**Theorem 1.11.** *The lattice of types* $(\mathcal{T}, \leq)$ *is anti–isomorphic to the lattice of rational cotorsion theories* $(\mathcal{C}_{rat}, \leq)$ *via the mapping*
$$\tau = t(T) \stackrel{\Phi}{\longmapsto} \left( {}^\perp(T^\perp), T^\perp \right).$$

With the above theorem we have fully described the lattice of all rational cotorsion theories.

Before we turn our attention to the general lattice of all cotorsion theories we need some "preparation", namely the construction of groups $G$ such that the cotorsion theories cogenerated by $G$, $\left( {}^\perp(G^\perp), G^\perp \right)$, are suitable for proving the Main Theorem. As the properties of these groups are interesting in their own right we consider them in a separate section.

## 2. An $\aleph_1$–free group without "small" cotorsion–free quotients

In this section we construct an $\aleph_1$–free group $G$ which has no proper subgroups $U$ of the same cardinality such that the quotient $G/U$ is cotorsion–free. In particular, if an epimorphic image $G/K$ of $G$ is cotorsion–free then the kernel $K$ is "small", namely $|K| < |G|$, and so $|G/K| = |G|$, i.e. the quotient is "big". Recall that an abelian group $G$ is said to be $\aleph_1$–free if all its countable subgroups are free or, equivalenty (by Pontryagin's criterion), if any finite rank subgroup is free. Note that $\aleph_1$–freeness implies cotorsion–freeness where a group is cotorsion–free if it is torsion–free and doesn't contain a copy of the rationals $\mathbb{Q}$ or the $p$–adic integers $J_p$ for some prime $p$.

In 1985 the first and second authors have constructed a cotorsion–free group $G$ with the above property where the cardinality of $G$ was a strong limit singular cardinal of cofinality bigger than $\aleph_0$ (see [11]). They already mentioned that using an at this stage newly developed set–theoretic method, which is nowadays known as Black Box, it is possible to replace the strong limit cardinal by any cardinal $\kappa$ with $\kappa^{\aleph_0} = \kappa$. This is what we will basically do here but to keep things simpler we only consider regular cardinals $\kappa$. Actually, we shall construct the group $G$ depending on a stationary set $S$ of $\kappa$ since that is what we need later in §3. However, the reader who is mainly interested in the construction of $G$ rather than in the application in §3 can ignore the statements regarding the stationary set. In particular, it would be enough to use the "ordinary" Black Box rather than the stationary one in order to obtain the required result. We need to use the stationary Black Box but we shall use an easy version of it.

Throughout this section let $\kappa$ be an infinite cardinal such that $|\alpha|^{\aleph_0} < \kappa$ for all ordinals $\alpha \in \kappa$ (e.g. take $\kappa = \mu^+$ for some $\mu$ with $\mu^{\aleph_0} = \mu$). Moreover, we fix a stationary set $S$ of $\kappa$ consisting of ordinals of cofinality $\omega$. First we define the parameters which are needed to formulate the Black Box.

Let $B$ be a free abelian group of rank $\kappa$, say $B = \bigoplus_{\alpha < \kappa} e_\alpha \mathbb{Z}$, and let $\widehat{B}$ denote the $p$–adic completion of $B$ for some fixed prime $p$. An element $b \in \widehat{B}$ can uniquely be written in the form $b = \sum_{\alpha < \kappa} e_\alpha b_\alpha$ where $b_\alpha \in J_p$. Thus we may define the *support* of $b \in \widehat{B}$ by $[b] = \{\alpha < \kappa \mid b_\alpha \neq 0\} \subseteq \kappa$; obviously $[b]$ is at most countable. This definition can be extended to subsets $M$ of $\widehat{B}$: $[M] = \bigcup_{b \in M} [b]$. Moreover, we define the *norm* of a subset $X$ of $\kappa$ by $\|X\| = \sup X = \sup_{x \in X} x$; this induces a norm for the elements $b$ and subsets $M$ of $\widehat{B}$: $\|b\| = \|[b]\|$, $\|M\| = \|[M]\| = \sup_{b \in M} \|b\|$.

While the Black Box, as for example known from [2], has been formulated using an "ordinary" tree $T = {}^{\omega >}\kappa$ and branches of this tree we need a different setting.

**Definition 2.1.** *A sequence $\overline{f} = (f_n)_{n \in \omega}$ of elements of $\widehat{B}$ is said to be the basis of a Signac–tree if*

(i) *$f_n$ is a pure element of $\widehat{B}$ for each $n \in \omega$,*
(ii) *$[f_n] \cap [f_m] = \emptyset$ for any $n \neq m$ in $\omega$,*
(iii) *$\|f_n\| < \|f_{n+1}\|$ for all $n \in \omega$.*

*Moreover, we call a subset $\{f_n \mid n \in X \subseteq \omega\}$ of $\overline{f}$ a branch over $\overline{f}$ and the set $\mathcal{S}(\overline{f})$ of all branches is said to be the Signac–tree over $\overline{f}$.*

Note that considering an "ordinary" tree $T = {}^{\omega >}\kappa$ the "basis of $T$" can be thought of as the set of all elements of $T$ of length (domain) 1 or, equivalently, as the elements of $\kappa$. Also note that the trees painted by the "pointilist" Signac are "dotted" which explains the name.

Next we define a relation on the set of all bases of Signac–trees:

**Definition 2.2.** *Let $\overline{f} = (f_n)_{n \in \omega}$ and $\overline{g} = (g_n)_{n \in \omega}$ be bases of Signac–trees. We say that $\overline{f}$ and $\overline{g}$ are* close *to each other if $\|\overline{f}\| = \|[\overline{f}] \cap [\overline{g}]\| = \|\overline{g}\|$. (Notation: $\overline{f} \sim \overline{g}$.)*

Note that the above defined relation is obviously an equivalence relation.

Since the Black Box is mainly a suitable enumeration of "traps" we need to say what we mean by it. Of course, the definition is adapted to our situation.

**Definition 2.3.** *A quadruple $\tau = (\overline{f}, P, K, b)$ is said to be a* trap *if*

(i) *$P$ is a canonical module, i.e. $P = \bigoplus_{\alpha \in X} e_\alpha \mathbb{Z}$ for some countable $X \subseteq \kappa$,*

  (ii) $\overline{f}$ is the basis of a Signac–tree of elements of $\widehat{P}$,
  (iii) $K$ is a countable pure subgroup of $\widehat{P}$,
  (iv) $b$ is a pure element of $P$, and
  (v) $\|b\| < \|P\| = \|K\| = \|\overline{f}\|$.

Moreover, we define the norm of $\tau$ by $\|\tau\| = \|P\|$ and we call a trap $\tau = (\overline{f}, P, K, b)$ as above an $S$–trap if $\|\tau\| \in S$.

We are now ready to present a suitable version of the Black Box. For a proof we refer to [2], [8] and [12]. Note that the Black Box is very robust under changes of its setting; the only real concern is the cardinality of the set of all objects in question. The choice of $\kappa$ and Definition 2.3 guarantee that all needed cardinalities are bounded by $\kappa$.

**The Black Box Lemma 2.4.** *There are an ordinal $\kappa^* < \kappa^+$ and a sequence $\left(\tau^\alpha = \left(\overline{f}^\alpha, P^\alpha, K^\alpha, b^\alpha\right)\right)_{\alpha<\kappa^*}$ of $S$–traps such that*

  (i) $\|\tau^\alpha\| \leq \|\tau^\beta\|$ *for* $\alpha < \beta$,
  (ii) $\left\|\left[\overline{f}^\alpha\right] \cap \left[\overline{f}^\beta\right]\right\| < \left\|\overline{f}^\beta\right\|$ *for* $\alpha < \beta$, *and*
  (iii) *for any pure submodule $U$ of $\widehat{B}$, any basis of a Signac–tree $\overline{g}$ of elements of $U$, and any pure element $b$ of $B$ with $\|b\| < \|\overline{g}\| \in S$ there is an $\alpha < \kappa^*$ such that*
$$K^\alpha \subseteq U \cap \widehat{P^\alpha}, \ g_n \in K^\alpha \ (n \in \omega), \ \overline{g} \sim \overline{f}^\alpha, \ \text{and} \ b = b^\alpha$$
  *where $(g_n)_{n\in\omega} = \overline{g}$.*

Next we construct the desired group $G$. We shall obtain $G$ lying between the free module $B$ and its completion $\widehat{B}$ by adding elements determined by suitable infinite branches of Signac–trees to $B$.

**Construction 2.5.** Let $\left(\tau^\alpha = \left(\overline{f}^\alpha, P^\alpha, K^\alpha, b^\alpha\right)\right)_{\alpha<\kappa^*}$ be a sequence of $S$–traps as in the Black Box Lemma 2.4. We construct $G = \bigcup_{\alpha<\kappa^*} G^\alpha$ inductively.

Let $G^0 = B$ and let $G^\alpha = \bigcup_{\beta<\alpha} G^\beta$ whenever $\alpha$ is a limit ordinal.

Now let $G^\alpha$ be given. If there is a basis of a Signac–tree $\overline{g}$ of elements of $K^\alpha$ which is close to $\overline{f}^\alpha$ then let $\overline{g}^\alpha = \overline{g}$; if this is not possible we put $\overline{g}^\alpha = \overline{f}^\alpha$. In the first case we call $\alpha$ a *strong ordinal* and in the latter we call it *weak ordinal*. In either case we define

$$G^{\alpha+1} = \left(G^\alpha + \sum_{\pi \in J_p} y_\pi^\alpha \mathbb{Z}\right)_*$$

where the index "$*$" denotes the purification within $\widehat{B}$ and the elements $y_\pi^\alpha$ are defined as follows: For each $p$–adic number $\pi \in J_p$ we choose an infinite branch $v_\pi = \{g_n^\alpha \mid n \in X_\pi\}$ ($X_\pi \subseteq \omega$ infinite) over $\overline{g}^\alpha = (g_n^\alpha)_{n\in\omega}$ such that $\omega \setminus (X_\pi \cap X_\rho)$ is infinite whenever $\pi \neq \rho$. Then let $a_\pi^\alpha = \sum_{n\in X_\pi} g_n^\alpha p^n$ and $y_\pi^\alpha = b^\alpha \pi + a_\pi^\alpha \in \widehat{B}$ ($\pi \in J_p$).

We can describe the above purification more explicitly:

For $\pi = \sum_{n\in\omega} a_n p^n \in J_p$ let $\pi_k = \sum_{n\geq k} a_n p^{n-k}$; let ${}_k a_\pi^\alpha = \sum_{n\in X_\pi,\, n\geq k} g_n^\alpha p^{n-k}$ and let ${}_k y_\pi^\alpha = b^\alpha \pi_k + {}_k a_\pi^\alpha$ ($k \in \omega$). Then we clearly have:

$$G_{\alpha+1} = G_\alpha + \sum_{\pi \in J_p,\, k \in \omega} {}_k y_\pi^\alpha \mathbb{Z}.$$

Finally, let $G = \bigcup_{\alpha < \kappa^*} G^\alpha$.

First note that $G$ is obviously a pure subgroup of $\widehat{B}$ of cardinality $\kappa = |\kappa^*|$.

We proceed with proving other properties of $G$. Next we show that $G$ is $\aleph_1$–free. In fact, we show more than that since the following proposition shall be used in §3. First we need:

**Lemma 2.6.** *Let $G$ be as constructed in 2.5.
Then $G_\alpha := \{g \in G \mid \|g\| < \alpha\}$ $(\alpha < \kappa)$ defines a $\kappa$–filtration of $G$.
Moreover, $G_{\alpha+1}/G_\alpha$ contains a non–zero $p$–divisible subgroup for $\alpha \in S$ and it is free otherwise.*

Note that we use the lower index $\alpha$ for the filtration $(\alpha \in \kappa)$ while we used the upper index $\alpha$ for the construction $(\alpha \in \kappa^*)$.

*of Lemma 2.6.* The first part of the result is obvious since $\|g\| < \alpha$ implies $g \in \widehat{\bigoplus_{\beta<\alpha} e_\beta \mathbb{Z}}$ and $\left|\widehat{\bigoplus_{\beta<\alpha} e_\beta \mathbb{Z}}\right| = |\alpha|^{\aleph_0} < \kappa$ for all $\alpha < \kappa$ by assumption.

So it remains to prove the second part. We consider the quotient group $G_{\alpha+1}/G_\alpha = \langle g + G_\alpha \mid g \in G,\ \|g\| = \alpha \rangle$.

If $\alpha \notin S$ then the elements of the form $g = x + e_\alpha z$ with $\|x\| < \alpha$, $z \in \mathbb{Z}$ are the only elements of norm $\alpha$. Hence $G_{\alpha+1}/G_\alpha = \langle e_\alpha + G_\alpha \rangle \cong \mathbb{Z}$ in this case.

For $\alpha \in S$ there is at least one $\gamma < \kappa^*$ such that $\alpha = \|\tau^\gamma\| = \|\overline{g}^\gamma\| = \|\overline{f}^\gamma\| = \|{}_k y_\pi^\gamma\|$ for $k \in \omega$, $\pi \in J_p$ since all elements of $S$ appear as norms in the sequence of $S$–traps $(\tau^\gamma)_{\gamma \in \kappa^*}$ (see 2.4 (iii)). But $\|g_n^\gamma\| < \|\overline{g}^\gamma\| = \alpha$ for all $n \in \omega$ and $\|b^\gamma\| < \alpha$ and so $y_\pi^\gamma - {}_k y_\pi^\gamma p^k = \sum_{n \in X_\pi,\, n<k} g_n^\gamma p^n + b^\gamma \sum_{n<k} a_n p^n \in G_\alpha$ for each $k \in \omega$ where $\sum_{n \in \omega} a_n p^n = \pi$. Therefore $y_\pi^\gamma$ is divisible by $p^k$ modulo $G_\alpha$ for each $k$ and thus $G_{\alpha+1}/G_\alpha$ contains a $p$–divisible subgroup for each $\alpha \in S$. □

**Proposition 2.7.** *Let $G = \bigcup_{\alpha<\kappa} G_\alpha$ be the $\kappa$–filtration of $G$ as in Lemma 2.6.
Then $G/G_\alpha$ is $\aleph_1$–free if and only if $\alpha \notin S$.*

*Proof.* By Lemma 2.6 we have that $G/G_\alpha \supseteq G_{\alpha+1}/G_\alpha$ contains a $p$–divisible subgroup whenever $\alpha \in S$ and thus it is not $\aleph_1$–free in this case.

So, let $\alpha \in \kappa \setminus S$. We show inductively that $G/G_\alpha = \bigcup_{\beta>\alpha} G_\beta/G_\alpha$ is $\aleph_1$–free.

Obviously, if $\beta$ is a limit ordinal and if $G_\gamma/G_\alpha$ is $\aleph_1$–free for each $\alpha < \gamma < \beta$ then $G_\beta/G_\alpha = \bigcup_{\alpha<\gamma<\beta} G_\gamma/G_\alpha$ is also $\aleph_1$–free by Pontryagin's criterion.

Moreover, $G_{\alpha+1}/G_\alpha$ is free and also $(G_{\beta+1}/G_\alpha)/(G_\beta/G_\alpha) \cong G_{\beta+1}/G_\beta$ is free for $\beta \notin S$ by Lemma 2.6 and hence $G_{\beta+1}/G_\alpha$ is $\aleph_1$–free provided $G_\beta/G_\alpha$ is $\aleph_1$–free.

Finally assume that $\beta \in S$ and that $G_\beta/G_\alpha$ is $\aleph_1$–free.
Let $X_\beta = \{\gamma < \kappa^* \mid \|\tau^\gamma\| = \beta\}$. Then $G_{\beta+1}/G_\alpha = (G_\beta/G_\alpha) + \langle {}_k y_\pi^\gamma + G_\alpha \mid k \in \omega,\ \pi \in J_p,\ \gamma \in X_\beta \rangle + \langle e_\beta + G_\alpha \rangle$.

We have to show that any finite set of elements of $G_{\beta+1}/G_\alpha$ is contained in a free pure subgroup of $G_{\beta+1}/G_\alpha$. Clearly, ${}_k y_\pi^\gamma \in \widehat{G_\beta}$ for all $\gamma \in X_\beta$, $\pi \in J_p$, $k \in \omega$, i.e. $[{}_k y_\pi^\gamma] \subseteq [G_\beta]$ and so we can ignore $e_\beta + G_\alpha$ since it is independent from all the other elements as $\alpha < \|e_\beta\| = \beta \notin [G_\beta]$.

Now, any finite subset of $G_{\beta+1}/G_\alpha$ is contained in a finite set $U = U_1 \cup U_2$ with $U_1 \subseteq G_\beta/G_\alpha$ and $U_2$ is of the form
$U_2 = \{{}_k y_\pi^\gamma + G_\alpha \mid k \leq l,\ \pi \in M,\ \gamma \in X\}$ where $l \in \omega$, and $M \subseteq J_p$, $X \subseteq X_\beta$ are finite sets.

Since $\beta = \lim_{n \in \omega} \|g_n^\gamma\|$ is a limit ordinal $(\gamma \in X_\beta)$ with $\alpha < \beta$ and $\|U_1\| < \beta$ we can find $l' \geq l$ such that $\alpha, \|U_1\| < \|g_{l'}^\gamma\|$ for all $\gamma \in X$. For each $\pi \in M$ let $l_\pi \geq l'$ be minimal with $l_\pi \in X_\pi$ and let
$U_1^* = \langle U_1, b^\gamma + G_\alpha, g_n^\gamma + G_\alpha \mid \gamma \in X, n \leq l' \rangle_* \subseteq G_\beta/G_\alpha$. Then $U$ is contained in the pure subgroup $U' = U_1^* \oplus \sum_{\pi \in M, \gamma \in X}({}_{l_\pi} y_\pi^\gamma + G_\alpha)\mathbb{Z}$ of $G_{\beta+1}/G_\alpha$ where $U_1^*$ is a pure finite rank subgroup of $G_\beta/G_\alpha$ and so $U_1^*$ is free by assumption.

It remains to show that $U_2^* := \sum_{\pi \in M,\ \gamma \in X}({}_{l_\pi} y_\pi^\gamma + G_\alpha)\mathbb{Z}$ is free. Since $\overline{g}^\gamma \sim \overline{f}^\gamma$ condition (ii) in the Black Box Lemma 2.4 implies $U_2^* =$
$\bigoplus_{\gamma \in X}\left(\sum_{\pi \in M}({}_{l_\pi} y_\pi^\gamma + G_\alpha)\mathbb{Z}\right)$ because starting with a maximal $\gamma_m \in X$ we can find $n_0 \in \omega$ such that $[g_n^{\gamma_m}] \not\subseteq \bigcup_{\gamma_m > \gamma \in X}[\overline{g}^\gamma]$ for all $n \geq n_0$ and we can proceed like this with the maximal element of the remaining set $X \setminus \{\gamma_m\}$ and so on.

Moreover, by our choice of the branches $v_\pi$ $(\pi \in J_p)$ in the Construction 2.5 we have that $\{{}_{l_\pi} y_\pi^\gamma + G_\alpha \mid \pi \in M\}$ $(\gamma \in X)$ is linearly independent since, for a fixed $\pi \in M$, we can recursively find integers $l_\pi \leq n_1 < n_2 < \ldots < n_k$ and sets $M \supset M_1 \supset M_2 \supset \ldots \supset M_k = \{\pi\}$ $(k \in \omega)$ such that $n_i \in X_\rho$ exactly if $\rho \in M_i$, i.e. the supports are sufficiently different. Hencefore $U_2^* \subseteq G_{\beta+1}/G_\alpha$ is free and so $U' = U_1^* \oplus U_2^*$ is free which completes the proof. $\square$

The most interesting property of the group $G$ is that it has no "small" cotorsion–free quotients. So $G$ as constructed in 2.5 is a suitable candidate for proving the final result of this section:

**Theorem 2.8.** *Let $\kappa$ be a regular cardinal with $|\alpha|^{\aleph_0} < \kappa$ for all $\alpha < \kappa$. Then there exists an $\aleph_1$–free group $G$ of cardinality $\kappa$ such that, for any subgroup $U$ of $G$ with $|U| = \kappa$, either $G = U$ or $G/U$ is not cotorsion–free.*

To prove this theorem we need a special case of the generalized $\Delta$–Lemma of Erdös and Rado [7]; for the proof see also [3].

**$\Delta$–Lemma 2.9.** *Let $\kappa$ be as in Theorem 2.8 and let $\Sigma$ be a family of cardinality $\kappa$ consisting of countable subsets of $\kappa$.
Then there is a subfamily $\Sigma'$ of $\Sigma$, also of cardinality $\kappa$, and an at most countable subset $F$ of $\kappa$ such that $X \cap Y = F$ for all $X \neq Y \in \Sigma'$.*

*of Theorem 2.8.* Let $B \subseteq_* G \subseteq_* \widehat{B}$ be as constructed in 2.5. Then $G$ is of cardinality $\kappa$. Moreover, $G = G/G_0$ is $\aleph_1$–free by Proposition 2.7.

Now let $U$ be a subgroup of $G$ of cardinality $\kappa$ such that $G/U$ is cotorsion–free. Moreover, let $\varphi: G \longrightarrow G/U$ be the canonical epimorphism.

If $B = \bigoplus_{\alpha < \kappa} e_\alpha \mathbb{Z} \subseteq U$ then $G/U$ is $p$–divisible as an epimorphic image of $G/B \subseteq_* \widehat{B}/B$. But $G/U$ is reduced as it is cotorsion–free and thus $G = U$ follows in this case.

Now assume that $\varphi$ is non zero, i.e. $U \neq G$. Then $B \not\subseteq U$, in particular there is a pure element $b \in B$ with $b\varphi \neq 0$ since $U \subseteq_* G$. We shall make use of this element $b$ later. First we apply the $\Delta$–Lemma 2.9 to the set $\Sigma = \{[u] \mid u \in U\}$. Then there is a subset $\Sigma'$ of $\Sigma$ of cardinality $\kappa$ and an at most countable set $F \subseteq \kappa$ such that $X \cap Y = F$ for all $X \neq Y \in \Sigma'$. For each $X \in \Sigma'$ choose one and only one element $u$ of $U$ with $X = [u]$ and let $U_1$ be the set of all such $u$'s. Then $|U_1| = \kappa$ and $[u] \cap [v] = F$ for all $u \neq v \in U_1$. Each $u \in U_1$ has a unique component $u \restriction F \in \widehat{\bigoplus_{\alpha \in F} e_\alpha \mathbb{Z}} \sqsubset \widehat{B}$. But $\left| \widehat{\bigoplus_{\alpha \in F} e_\alpha \mathbb{Z}} \right| = |F|^{\aleph_0} = 2^{\aleph_0} < \kappa = |U_1|$ and thus there are a subset $U_2$ of $U_1$ also of cardinality $\kappa$ and an element $f$ of $\widehat{\bigoplus_{\alpha \in F} e_\alpha \mathbb{Z}}$ such that $u \restriction F = f$ for all $u \in U_2$. Next let $U_3$ and $U_4$ be disjoint subsets of $U_2$ with $U_2 = U_3 \cup U_4$, both of cardinality $\kappa$, and let $U_5 = \{x - y \mid x \in U_3, y \in U_4\}$. Then $|U_5| = \kappa$ and $[u] \cap [v] = \emptyset$ for all $u \neq v \in U_5$. Replace each element $u \in U_5$ by its "purification", i.e. by $up^{-n_u}$ where $n_u$ is the maximal power of $p$ dividing $u$; clearly $[u] = [up^{-n_u}]$. Call the new set $U_6$.

Now let $U^* = \{u \in U_6 \mid \|u\| > \|b\|\}$ where $b$ is the pure element from above with $b\varphi \neq 0$. Then $U^*$ consists of pure elements of $U = \ker \varphi$ and satisfies $[u] \cap [v] = \emptyset$ for all $u \neq v \in U^*$. Moreover, $|U^*| = \kappa$ since the set of all elements of $\widehat{B}$ with norm less than or equal to $\|b\| = \beta$ is of cardinality $|\beta|^{\aleph_0}$.

We consider the set $\mathcal{S}$ of all sequences $\overline{g} = (g_n)_{n \in \omega}$ of elements of $U^*$ with $\|g_n\| < \|g_{n+1}\|$ for all $n \in \omega$. Then $|\mathcal{S}| = \kappa$ and the set $C = \left\{ \|\overline{g}\| = \sup_{n \to \omega} \|g_n\| \mid \overline{g} \in \mathcal{S} \right\}$ is unbounded and closed under limits of countable subsets; a set $C$ satisfying these properties is called an $\omega$–cub. Since our fixed stationary set $S$ consists of ordinals of cofinality $\omega$ the intersection with any $\omega$–cub is non empty: $C \cap S \neq \emptyset$. Therefore there is an element $\overline{g} = (g_n)_{n \in \omega}$ of $\mathcal{S}$ with $\|\overline{g}\| \in S$. Obviously, $\overline{g}$ is a basis of a Signac–tree since $g_n \in U^*$ and $\|g_n\| < \|g_{n+1}\|$ for all $n \in \omega$ (see Definition 2.1). By the Black Box Lemma 2.4 there is an ordinal $\alpha < \kappa^*$ such that
$K^\alpha \subseteq U \cap P^\alpha$, $g_n \in K^\alpha$ $(n \in \omega)$, $\overline{g} \sim \overline{f}^\alpha$ and $b = b^\alpha$. Therefore $\alpha$ is a strong ordinal.

We now consider the elements $y_\pi^\alpha = b\pi + a_\pi^\alpha$ $(\pi \in J_p)$. Since $g_n^\alpha \in K^\alpha \subseteq U = \ker \varphi$ $(n \in \omega)$ and $a_\pi^\alpha = \sum_{n \in X_\pi} g_n^\alpha p^n$ $(\pi \in J_p)$ we have $a_\pi^\alpha \widehat{\varphi} = 0$ where the continuous homomorphism $\widehat{\varphi} : \widehat{G} = \widehat{B} \longrightarrow \widehat{(G/U)}$ is the unique extension of $\varphi$. The continuity of $\widehat{\varphi}$ also implies $y_\pi^\alpha \varphi = (b\pi + a_\pi^\alpha)\varphi = (b\pi)\widehat{\varphi} = b\widehat{\varphi}\pi = b\varphi\pi \in G/U$ for all $\pi \in J_p$ and so $(b\varphi)J_p \subseteq G/U$, i.e. $G/U$ contains a copy of the $p$–adic integers $J_p$ contradicting the cotorsion–freeness of $G/U$. This implies $G = U$ and so the proof is finished. $\square$

Note that additional to the above properties of $G$ we could prescribe the endomorphism ring of $G$ using the standard methods via the Black Box.

However, the group $G$ as constructed in this section is exactly what we need to prove the Main Theorem in §3. In fact, we shall need a family of such groups $G$ depending on different stationary sets $S$; thus we introduce the notation $G = G(S)$ to refer to the group $G$ as constructed in 2.5 satisfying the conclusion of Theorem 2.8.

3. Embeddings of posets into the lattice of cotorsion theories

Throughout this section let $I$ be an arbitrary set and let $\mathcal{P} = \mathcal{P}(I)$ be the power set of $I$. Moreover, let $\kappa \geq |I|$ be a regular cardinal such that, for all ordinals $\alpha < \kappa$, $|\alpha|^{\aleph_0} < \kappa$. Note that such a cardinal always exists, e.g. take $\kappa = \left(|I|^{\aleph_0}\right)^+$.
The aim of this section is to prove the Main Theorem of the paper:

**Main Theorem 3.1.** *There is an embedding from $(\mathcal{P}, \subseteq)$ into the lattice of all cotorsion theories $(\mathcal{C}, \leq)$.*

Note that any poset can be embedded into the power set lattice of some set $I$.

We shall prove the main theorem in several steps. First we define an order–reversing mapping $\Phi : (\mathcal{P}, \subseteq) \longrightarrow (\mathcal{C}, \leq)$ which will turn out to be injective. Since the mapping $(\mathcal{P}, \subseteq) \longrightarrow (\mathcal{P}, \subseteq)$ $(X \mapsto I \setminus X)$ is an order–reversing isomorphism this induces the required embedding.

Now, the set $\{\alpha \in \kappa \mid cf(\alpha) = \omega\}$ is stationary and can be partitioned into $|I|$ disjoint stationary subsets; say $\{\alpha \in \kappa \mid cf(\alpha) = \omega\} = \bigcup_{i \in I} S_i$. Let $G_i = G(S_i)$ be an $\aleph_1$–free group of cardinality $\kappa$ as constructed in 2.5 depending on the stationary set $S_i$ $(i \in I)$. Moreover, for each $X \subseteq I$, let $G_X = \bigoplus_{i \in X} G_i$. We define $X\Phi = (^\perp(G_X^\perp), G_X^\perp) \in \mathcal{C}$. Obviously, $\Phi$ is well defined and, for $Y \subseteq X \subseteq I$, we have $G_Y \subseteq G_X$ and thus $G_Y^\perp \supseteq G_X^\perp$, i.e. $\Phi$ is order reversing. Note that $G_\emptyset = 0$ and so $\emptyset\Phi = (\mathcal{L}, \text{Mod–}\mathbb{Z})$ is the maximal cotorsion theory; recall that $\mathcal{L}$ denotes the class of all free abelian groups.

In order to establish that $\Phi$ is injective we construct groups $H^X$ ($\emptyset \neq X \subseteq I$) such that $\text{Ext}(G_Y, H^X) = 0$ if and only if $Y \subseteq X$, i.e. if $Y \not\subseteq X$ then $H^Y \in G_Y^\perp \setminus G_X^\perp$. The construction is based on a method of making Ext vanish which has been introduced by the first and second authors [14]; here we use the generalized method as developed by Eklof and Trlifaj [6].

**Construction 3.2.** Let $X$ be a fixed non–empty subset of $I$ and let $\lambda$ be a cardinal with $\lambda^\kappa = \lambda$. Moreover, let $H$ be a set of cardinality $\lambda$ with a $\lambda$–filtration $H = \bigcup_{\alpha < \lambda} H_\alpha$ such that $|H_0| = \kappa$ and $|H_\alpha| = |\alpha| \cdot \kappa = |H_{\alpha+1} \setminus H_\alpha|$ for all $\alpha < \lambda$.
We inductively define a group structure on $H$ and call the obtained group $H^X$.

We fix free resolutions $0 \longrightarrow K_i \longrightarrow F_i \longrightarrow G_i \longrightarrow 0$ of $G_i$ with $|K_i| = |F_i| = \kappa$ $(i \in X)$ and we "enumerate" all set mappings from all $K_i$'s into $H$ by $\bigcup_{i \in X} {}^{K_i}H = \{\varphi_\alpha \mid \alpha < \lambda\}$ in such a way that each mapping appears $\lambda$ times.

Now let $H_0^X = \mathbb{Z}^{(\kappa)}$ be a free group of rank $\kappa$. If $\alpha$ is a limit ordinal and if the group structure $H_\beta^X$ on $H_\beta$ is defined for all $\beta < \alpha$ such that $H_\beta^X$ is a subgroup of $H_{\beta+1}^X$ then let $H_\alpha^X = \bigcup_{\beta < \alpha} H_\beta^X$ have the induced group structure.

Now let the group structure $H_\alpha^X$ be given.

If $\text{Im}\varphi_\alpha \subseteq H_\alpha^X$ and if $\varphi_\alpha$ is a homomorphism then let $\widetilde{\varphi_\alpha} = \varphi_\alpha$, and put $\widetilde{\varphi_\alpha} = 0$ otherwise. In either case we define $H_{\alpha+1}^X$ to be the pushout

$$\begin{array}{ccc} F_i & \xrightarrow{\psi_\alpha} & H_{\alpha+1}^X \\ \uparrow & & \uparrow \\ K_i & \xrightarrow{\widetilde{\varphi_\alpha}} & H_\alpha^X \end{array}$$

where $dom\varphi_\alpha = K_i$ for some $i \in X$. Hence $\psi_\alpha$ is an extension of $\widetilde{\varphi_\alpha}$ and $H^X_{\alpha+1}/H^X_\alpha \cong F_i/K_i \cong G_i$.
Finally let the structure on $H^X = \bigcup_{\alpha < \lambda} H^X_\alpha$ be the induced one.

Note that the cardinality of $H^X$ is obviously $\lambda$ for each non–empty set $X \subseteq I$.
First we show that $H^X \in G_Y^\perp$ for any set $Y \subseteq X$.

**Proposition 3.3.** Let $\emptyset \neq X \subseteq I$ and let $H^X$ be as constructed in 3.2.
Then $H^X \in G_Y^\perp$, i.e. $\mathrm{Ext}(G_Y, H^X) = 0$, for any $Y \subseteq X$.
*Proof.* Since $G_Y = \bigoplus_{i \in Y} G_i$ $(Y \subseteq X)$ it is sufficient to show that
$\mathrm{Ext}(G_i, H^X) = 0$ for each $i \in X$.
We consider the free resolution $0 \longrightarrow K_i \longrightarrow F_i \longrightarrow G_i \longrightarrow 0$ of $G_i$ as in Construction 3.2. Let $\varphi : K_i \longrightarrow H^X = \bigcup_{\alpha < \lambda} H^X_\alpha$ be a homomorphism. Since $|K_i\varphi| \leq |K_i| = \kappa < cf(\lambda)$ there is an ordinal $\beta < \lambda$ such that $\mathrm{Im}\varphi \subseteq H^X_\beta$. Moreover, by the enumeration of $\bigcup_{i \in X} {}^{K_i}H$ in 3.2 there is $\beta \leq \alpha < \lambda$ such that $\varphi = \varphi_\alpha = \widetilde{\varphi_\alpha}$ and thus there is an extension $\psi_\alpha : F_i \longrightarrow H^X$ of $\varphi$.
Therefore we have seen that every homomorphism from $K_i$ into $H^X$ extends to a homomorphism from $F_i$ into $H^X$ and hence $\mathrm{Ext}(G_i, H^X) = 0$ whenever $i \in X$. This implies $\mathrm{Ext}(G_Y, H^X) = 0$ for all $Y \subseteq X$. $\square$

It remains to show that $\mathrm{Ext}(G_i, H^X) \neq 0$ whenever $i \notin X$. Although it seems to be the more likely case that $\mathrm{Ext}(A, B) \neq 0$ for arbitrary groups $A$ and $B$ there is some work to do in order to prove this. The key to the proof is the following:

**Definition 3.4.** Let $\kappa$ be as above and let $S$ be a stationary set in $\kappa$.
We call a group $A$ *locally S–free* if, for any smooth ascending chain $\{K_\alpha \mid \alpha < \kappa\}$ of subgroups $K_\alpha$ of $A$ with $|K_\alpha| < \kappa$ for all $\alpha < \kappa$, the set $\{\delta \in S \mid K_{\delta+1}/K_\delta \text{ not } \aleph_1\text{–free}\}$ is not stationary in $\kappa$.

First we investigate the groups $G_i = G(S_i)$ $(i \in I)$ with respect to the just defined property.

**Proposition 3.5.** Let $i \neq j$ be elements of $I$ and let $G_i = G(S_i)$ be as above.
Then $G_i$ is locally $S_j$–free but not locally $S_i$–free.
*Proof.* By Lemma 2.6 and Proposition 2.7 there is a $\kappa$–filtration $G_i = \bigcup_{\alpha < \kappa} G^i_\alpha$ of $G_i$ such that

$$G_i/G^i_\alpha \text{ is } \aleph_1\text{–free} \quad \text{iff} \quad \alpha \notin S_i. \qquad (*)$$

Moreover, we know that $G^i_{\alpha+1}/G^i_\alpha$ contains a divisible subgroup for any $\alpha \in S_i$. Hence it follows immediately that $G_i$ is not locally $S_i$–free since $\{\delta \in S_i \mid G^i_{\delta+1}/G^i_\delta \text{ not } \aleph_1\text{–free}\} = S_i$ is stationary in $\kappa$.
Thus, it remains to show that $G_i$ is locally $S_j$–free for any $j \neq i$.
Let $\{K_\alpha \mid \alpha < \kappa\}$ be a smooth ascending chain of subgroups of $G_i$ all of cardinality less than $\kappa$ and let $K = \bigcup_{\alpha < \kappa} K_\alpha$.
If $|K| < \kappa$ then there is an $\alpha_0 < \kappa$ such that $K = K_\alpha = K_{\alpha_0}$ for all $\alpha > \alpha_0$. Thus the set $\{\delta \in S_j \mid K_{\delta+1}/K_\delta \text{ not } \aleph_1\text{–free}\}$ is bounded by $\alpha_0$ and hence it is not stationary.

Otherwise $|K| = \kappa$ and $\{K_\alpha \mid \alpha < \kappa\}$ is a $\kappa$–filtration of $K$. Also $\{K \cap G_\alpha^i \mid \alpha < \kappa\}$ is a $\kappa$–filtration of $K$ and thus there exists a closed unbounded set (cub) $C$ in $\kappa$ such that $K_\alpha = K \cap G_\alpha^i$ for all $\alpha \in C$. Let $\delta \in C \cap S_j$, then $K_\delta = K \cap G_\delta^i$ since $\delta \in C$ and $G_i/G_\delta^i$ is $\aleph_1$–free by $(*)$ since $\delta \in S_j$ which is disjoint from $S_i$. Therefore $K_{\delta+1}/K_\delta \subseteq K/K_\delta = K/(K \cap G_\delta^i) \cong (K + G_\delta^i)/G_\delta^i \subseteq G_i/G_\delta^i$ is $\aleph_1$–free. Hence $C$ is disjoint from $\{\delta \in S_j \mid K_{\delta+1}/K_\delta \text{ not } \aleph_1\text{–free}\}$ and thus this set is not stationary. So we have shown that $G_i$ is locally $S_j$–free. $\square$

As an immediate consequence from the above proposition we have:

**Corollary 3.6.** *Let $A$ be a cotorsion–free locally $S_i$–free group for some $i \in I$. Then $\mathrm{Hom}(G_i, A) = 0$.*
*Proof.* Suppose, for contradiction, that there exists a non–zero homomorphism $\varphi : G_i = \bigcup_{\alpha < \kappa} G_\alpha^i \longrightarrow A$ and let $K_\alpha = G_\alpha^i \varphi$. Since $G_i/\ker \varphi \cong \mathrm{Im}\,\varphi \subseteq A$ is cotorsion–free by assumption it follows from Theorem 2.8 that the kernel of $\varphi$ has to be "small", i.e. $|\ker \varphi| < \kappa$. Therefore there is $\alpha_0 < \kappa$ such that $\ker \varphi \subseteq G_\alpha^i$ for all $\alpha \geq \alpha_0$ as $\kappa$ is regular. This implies $K_{\alpha+1}/K_\alpha \cong (G_{\alpha+1}^i/\ker \varphi)/(G_\alpha^i/\ker \varphi) \cong G_{\alpha+1}^i/G_\alpha^i$ for all $\alpha \geq \alpha_0$ and thus, by Lemma 2.6, $\{\delta \in S_i \mid K_{\delta+1}/K_\delta \text{ not } \aleph_1\text{–free}\} = \{\delta \in S_i \mid \delta \geq \alpha_0\}$ is stationary contradicting the local $S_i$–freeness of $A$. $\square$

We now proceed with investigating the relevant properties of the groups $H^X$ ($\emptyset \neq X \subseteq I$) as constructed in 3.2. Since $G_i = G(S_i)$ is $\aleph_1$–free for all $i \in I$ by Theorem 2.8 we immediately have:

**Lemma 3.7.** *Let $\emptyset \neq X \subseteq I$ and let $H^X = \bigcup_{\alpha < \lambda} H_\alpha^X$ be as constructed in 3.2. Then $H^X$ and $H^X/H_\alpha^X$ ($\alpha < \lambda$) are $\aleph_1$–free.*

Next we consider the local $S_i$–freeness of $H^X$.

**Proposition 3.8.** *Let $\emptyset \neq X \subseteq I$, $i \in I \setminus X$ and let $H^X$ be as constructed in 3.2. Then $H^X$ and $H^X/H_0^X$ are locally $S_i$–free.*
*Proof.* We shall show inductively that $H_\alpha^X$ and $H_\alpha^X/H_0^X$ are locally $S_i$–free for all $\alpha < \lambda$. Of course, $H_0^X = \mathbb{Z}^{(\kappa)}$ and $H_0^X/H_0^X = 0$ are locally $S_i$–free. In the following we restrict our attention to the $H_\alpha^X$'s as there is no difference in the arguments when considering the $(H_\alpha^X/H_0^X)$'s.
First assume that $H_\alpha^X$ is locally $S_i$–free and consider a smooth ascending chain $\{K_\gamma \mid \gamma < \kappa\}$ of subgroups of $H_{\alpha+1}^X$ with $|K_\gamma| < \kappa$. Let $K = \bigcup_{\gamma < \kappa} K_\gamma$ and let
$M = \{\delta \in S_i \mid K_{\delta+1}/K_\delta \text{ not } \aleph_1\text{–free}\}$. Moreover, let
$M_1 = \{\delta \in S_i \mid (K_{\delta+1} \cap H_\alpha^X)/(K_\delta \cap H_\alpha^X) \text{ not } \aleph_1\text{–free}\}$ and
$M_2 = \{\delta \in S_i \mid (K_{\delta+1} + H_\alpha^X)/(K_\delta + H_\alpha^X) \text{ not } \aleph_1\text{–free}\}$.
By induction hypothesis, $M_1$ is not stationary. Also $M_2$ is not stationary by Proposition 3.5 since $\{(K_\gamma + H_\alpha^X)/H_\alpha^X \mid \gamma < \kappa\}$ is a chain in $H_{\alpha+1}^X/H_\alpha^X \cong G_j$ for some $j \in X$. Thus there are cubs $C_1$ and $C_2$ such that $M_l \cap C_l = \emptyset$ ($l = 1, 2$). Let $C = C_1 \cap C_2$, then $C$ is also a cub. We prove that $M \cap C = \emptyset$ and therefore $M$ is not stationary.
Let $\delta \in C \cap S_i$. Then $\delta \notin M_1 \cup M_2$ and thus $(K_{\delta+1} \cap H_\alpha^X)/(K_\delta \cap H_\alpha^X)$ and $(K_{\delta+1} + H_\alpha^X)/(K_\delta + H_\alpha^X)$ are $\aleph_1$–free. There is an epimorphism $K_{\delta+1}/K_\delta \longrightarrow (K_{\delta+1} + H_\alpha^X)/(K_\delta + H_\alpha^X)$ with kernel $((K_{\delta+1} \cap H_\alpha^X) +$ $K_\delta)/K_\delta \cong (K_{\delta+1} \cap H_\alpha^X)/(K_\delta \cap H_\alpha^X)$ and hence $K_{\delta+1}/K_\delta$ is $\aleph_1$–free as an extension of an $\aleph_1$–free group by an $\aleph_1$–free group. Therefore $\delta \notin M$ and so $M \cap C = \emptyset$.

Now let $\alpha$ be a limit ordinal and suppose that $H_\beta^X$ is locally $S_i$–free for all $\beta < \alpha$. Consider $K = \bigcup_{\gamma<\kappa} K_\gamma \subseteq H_\alpha^X$ with $|K_\gamma| < \kappa$ for all $\gamma < \kappa$. Moreover, let $M = \{\delta \in S_i \mid K_{\delta+1}/K_\delta \text{ not } \aleph_1\text{–free}\}$.

If $K \subseteq H_\beta^X$ for some $\beta < \alpha$ then $M$ is not stationary by assumption.

So assume otherwise, i.e. $K \nsubseteq H_\beta^X$ for all $\beta < \alpha$. Then the cofinality of $\alpha$ is less than or equal to $\kappa$.

First we consider the case of $cf(\alpha) = \mu < \kappa$. Let $\alpha = \sup_{\nu \to \mu} \alpha_\nu$ with $\alpha_\nu < \alpha$ and put $H'_\nu = H_{\alpha_\nu}^X$ ($\nu < \mu$). Then $H_\alpha^X = \bigcup_{\nu<\mu} H'_\nu$. By assumption there are cubs $C_\nu$ ($\nu < \mu$) such that $C_\nu \cap M_\nu = \emptyset$ where
$M_\nu = \{\delta \in S_i \mid (K_{\delta+1} \cap H'_\nu)/(K_\delta \cap H'_\nu) \text{ not } \aleph_1\text{–free}\}$. Now let $C = \bigcap_{\nu<\mu} C_\nu$. Then $C$ is also a cub by [5, II.4.3]. We show that $M \cap C = \emptyset$. Let $\delta \in C \cap S_i$. Then $(K_{\delta+1} \cap H'_\nu)/(K_\delta \cap H'_\nu) \cong ((K_{\delta+1} \cap H'_\nu) + K_\delta)/K_\delta$ is $\aleph_1$–free for each $\nu < \mu$ and so $K_{\delta+1}/K_\delta = \left(\bigcup_{\nu<\mu}(K_{\delta+1} \cap H'_\nu)\right)/K_\delta = \bigcup_{\nu<\mu}((K_{\delta+1} \cap H'_\nu) + K_\delta)/K_\delta$ is also $\aleph_1$–free. Hence $\delta \notin M$ and so
$M \cap C = \emptyset$ as required.

It remains to consider the case $cf(\alpha) = \kappa$. Since $|K_\gamma| < \kappa$ for each $\gamma < \kappa$ we may choose $\alpha_\gamma < \alpha$ such that $K_\gamma \subseteq H_{\alpha_\gamma}^X$ and $\alpha_\gamma = \sup_{\delta \to \gamma} \alpha_\delta$ whenever $\gamma$ is a limit. Then $\alpha = \sup_{\gamma \to \kappa} \alpha_\gamma$ since $K \nsubseteq H_\beta^X$ for any $\beta < \alpha$.

Now let $H'_\gamma = H_{\alpha_\gamma}^X$ and let
$M_\gamma = \{\delta \in S_i \mid (K_{\delta+1} \cap H'_\gamma)/(K_\delta \cap H'_\gamma) \text{ not } \aleph_1\text{–free}\}$ ($\gamma < \kappa$).
By assumption there are cubs $C_\gamma$ with $M_\gamma \cap C_\gamma = \emptyset$. We define $C$ to be the diagonal intersection $C = \Delta\{C_\gamma \mid \gamma < \kappa\} = \{\delta < \kappa \mid \delta \in \bigcap_{\gamma<\delta} C_\gamma\}$.

By [5, II.4.10] $C$ is also a cub. As before we show that $M \cap C = \emptyset$ in order to establish that $M$ is not stationary.

Consider an ordinal $\delta \in C \cap S_i$. Then $\delta$ is a limit ordinal and $\delta \in C_\gamma$ for all $\gamma < \delta$. So $(K_{\delta+1} \cap H'_\gamma)/(K_\delta \cap H'_\gamma) \cong ((K_{\delta+1} \cap H'_\gamma) + K_\delta)/K_\delta$ is $\aleph_1$–free ($\gamma < \delta$). It follows immediately that $(K_{\delta+1} \cap H'_\delta)/K_\delta = \bigcup_{\gamma<\delta}(((K_{\delta+1} \cap H'_\gamma) + K_\delta)/K_\delta)$ is also $\aleph_1$–free. Moreover, for all $\gamma \geq \delta$ we have
$((K_{\delta+1} \cap H'_{\gamma+1})/K_\delta)/((K_{\delta+1} \cap H'_\gamma)/K_\delta) \cong (K_{\delta+1} \cap H'_{\gamma+1})/(K_{\delta+1} \cap H'_\gamma) \cong ((K_{\delta+1} \cap H'_{\gamma+1}) + H'_\gamma)/H'_\gamma \subseteq H^X/H'_\gamma$ is $\aleph_1$–free. Therefore it follows by transfinite induction that $K_{\delta+1}/K_\delta = \bigcup_{\delta \leq \gamma < \kappa}((K_{\delta+1} \cap H'_\gamma)/K_\delta)$ is $\aleph_1$–free. Thus $\delta$ cannot be an element of $M$ which completes the proof. □

Using the above result we can finally prove the last missing bit in order to establish the correctness of the Main Theorem.

**Proposition 3.9.** *Let $\emptyset \neq X \subseteq I$, $i \in I \setminus X$, and let $H^X$ be as constructed in 3.2. Then $\mathrm{Ext}(G_i, H^X) \neq 0$, i.e. $H^X \notin G_i^\perp$.*

*Proof.* Let $0 \longrightarrow K_i \longrightarrow F_i \longrightarrow G_i \longrightarrow 0$ be a free resolution of $G_i$ with $|K_i| = |F_i| = \kappa$. In order to show $\mathrm{Ext}(G_i, H^X) \neq 0$ it is enough to find a homomorphism $\varphi : K_i \longrightarrow H^X$ which doesn't extend to a homomorphism $\widetilde{\varphi} : F_i \longrightarrow H^X$.

Let $\varphi : K_i \longrightarrow H_0^X \subseteq H^X$ be an isomorphism between the two free groups $K_i, H_0^X$ of rank $\kappa$. Suppose, for contradiction, that there is $\widetilde{\varphi} : F_i \longrightarrow H^X$ with $\widetilde{\varphi} \restriction K_i = \varphi$. Then $\widetilde{\varphi}$ induces a homomorphism $\overline{\varphi} : F_i / K_i \cong G_i \longrightarrow H^X / H_0^X$. But $H^X / H_0^X$ is $\aleph_1$–free by Lemma 3.7; in particular it is cotorsion–free. Also $H^X / H_0^X$ is locally $S_i$–free by Proposition 3.8. Hence $\mathrm{Hom}(G_i, H^X / H_0^X) = 0$ by Corollary 3.6. Therefore $\overline{\varphi} = 0$ and so $F_i \widetilde{\varphi} = H_0$. But this implies $F_i = K_i \oplus \ker \widetilde{\varphi}$ since, for each $f \in F_i$, there is $k \in K_i$ with $f\widetilde{\varphi} = k\widetilde{\varphi} = k\varphi$ and $\ker \widetilde{\varphi} \cap K_i = \ker \varphi = \{0\}$. Then $G_i \cong F_i / K_i \cong \ker \widetilde{\varphi}$ is free contradicting that $G_i$ is not free. Therefore there is no such extension $\widetilde{\varphi}$ of $\varphi$ and thus $\mathrm{Ext}(G_i, H^X) \neq 0$ as required. □

In the above results 3.2 – 3.9 we have shown that the mapping
$\Phi : (\mathcal{P}, \subseteq) \longrightarrow (\mathcal{C}, \leq)$ as defined at the beginning of this section is an order–reversing injection. Therefore we have proven the Main Theorem. Finally note that it follows immediately from the Main Theorem that there are ascending, descending, and anti–chains of arbitrary size in the lattice of all cotorsion theories and this answers the original question which led to this paper.

Fachbereich 6, Mathematik und Informatik,, Universität Essen,, 45117 Essen, Germany
*E-mail address*: `R.Goebel@Uni-Essen.De`

Department of Mathematics,, Hebrew University,, Jerusalem, Israel
*E-mail address*: `shelah@math.huji.ac.il`

Fachbereich 6, Mathematik und Informatik,, Universität Essen,, 45117 Essen, Germany
*E-mail address*: `Simone.Wallutis@Uni-Essen.De`